\documentclass[final]{siamltex}
\usepackage{amssymb,amsmath,amsfonts,color}
\usepackage{graphicx}
\usepackage{epstopdf}

\newcommand{\cblue}[1]{{\color{black}{#1}}}

\newcommand{\sxy}[2]{\left\langle #1,#2 \right\rangle}
\newtheorem{example}[theorem]{Example}

\newcommand{\piF}[1]{\pi_k}

\title{The field of values bounds on ideal GMRES}
\author{J\"org Liesen\footnotemark[1] and Petr Tich\'y\footnotemark[3]}

\begin{document}

\footnotetext[1]{Institute of Mathematics, Technical University of Berlin,
Stra{\ss}e des 17. Juni 136, 10623 Berlin, Germany ({\tt liesen@math.tu-berlin.de}).}

\footnotetext[3]{Faculty of Mathematics and Physics, Charles University, Sokolovsk\'a 83, 182 07 Prague 8, Czech Republic, ({\tt ptichy@karlin.mff.cuni.cz}).
This work was partially supported by project 17-04150J of the Grant Agency of the Czech Republic.}

\maketitle

%\begin{center}
%{\em To Zden\v{e}k Strako\v{s} on his 60th birthday}
%\end{center}\medskip

\begin{abstract}
A widely known result of Elman, and its improvements due to Starke, Eiermann and Ernst,
gives a bound on the worst-case GMRES residual norm using quantities related to the field
of values of the given matrix and its inverse. We prove that these bounds
also hold for the ideal GMRES approximation, and we derive \cblue{and discuss} some improvements of the bounds.
\end{abstract}



\pagestyle{myheadings}
\thispagestyle{plain}
\markboth{J. LIESEN and P. TICH\'Y}{FIELD OF VALUES BOUNDS ON IDEAL GMRES}

\section{Introduction}

Consider a linear algebraic system $Ax=b$ with a nonsingular matrix \cblue{$A\in\mathbb{C}^{n\times n}$ and a right hand side $b\in\mathbb{C}^n$.
Given an initial approximation $x_0\in \mathbb{C}^n$} and the initial residual $r_0\equiv b-Ax_0$,
the GMRES method of Saad and Schultz~\cite{SaSc1986} iteratively constructs approximations~$x_k$
such that
\begin{equation}\label{eqn:GMRES}
\|r_k\|=\|b-Ax_k\|=\min_{p\in \piF{F}}\|p(A)r_0\|,\quad k=1,2,\dots,
\end{equation}
where $\|v\|\equiv \langle v,v\rangle^{1/2}$ denotes the Euclidean norm on ${\mathbb C}^n$,
and $\piF{F}$ is the set of polynomials $p$ of degree at most $k$ with coefficients
in $\mathbb{C}$, and with $p(0)=1$. \cblue{For real data $A$ and $b$,
the coefficients of the polynomial $p$ solving \eqref{eqn:GMRES} are real, and, therefore, we can consider $\piF{F}$ to be the set of polynomials with real coefficients in this case.}

The convergence analysis of GMRES has been a challenge since the introduction of the algorithm;
see~\cite{LiTi2005} or~\cite[Section~5.7]{B:LiSt2013} for surveys of this research area. Here
we focus on GMRES convergence bounds that are independent of the initial residual, i.e., for a given $A$,
we consider
the worst-case behavior of the method. It is easy to see that for each given $A$, $b$ and $x_0$,
the $k$th relative GMRES residual norm satisfies
\begin{equation}\label{eqn:WC}
\frac{\|r_k\|}{\|r_0\|}\leq
\max_{v\in {\mathbb C}^n \atop \|v\|=1} \min_{p\in \piF{F}} \|p(A)v\|.
\end{equation}
{The expression on the right hand side is called the $k$th {\em worst-case GMRES
residual norm}. For each given matrix $A$ and iteration step $k$, this quantity is attainable
by the relative GMRES residual norm for some initial residual $r_0$.}
Mathematical properties of worst-case GMRES have been studied in~\cite{FaLiTi2013}; see
also~\cite{LiTi2004}.

\cblue{Consider now real data and} let $M\equiv\frac12 (A+A^T)$ be the symmetric part of $A$.
Assuming that $M$ is positive definite, a widely known result of Elman, stated originally
for the relative residual norm of the GCR method in~\cite[Theorem~5.4 and~5.9]{T:El1982},
implies that
\begin{equation}\label{eqn:Elman}
\max_{v\in {\mathbb R}^n \atop \|v\|=1} \min_{p\in \piF{R}} {\|p(A)v\|} \leq
\left(1-\frac{\lambda_{\min}(M)^2}{\lambda_{\max}(A^TA)}\right)^{k/2};
\end{equation}
see also the paper~\cite[Theorem~3.3]{EiElSc1983}.
%More recently, Elman's
%bound has been generalized in \cite{BeGoTy2005} for nonsingular complex matrices, and
%some asymptotically sharper bounds have been proposed. These will be discussed in more
%detail in Section~\ref{sec:improve} below.

Let ${\cal F}(A)$ be the field of values \cblue{of~$A$}, and
let $\nu(A)$ be the distance of ${\cal F}(A)$ from the origin, i.e.,
\begin{align*}
{\cal F}(A) &\equiv \{\left\langle Av,v\right\rangle : v\in\mathbb{C}^{n},\ \|v\|=1\}, \quad
\nu(A) \equiv \min_{z\in {\cal F}(A)}\,|z|.
\end{align*}
Then the bound \eqref{eqn:Elman} can be written as
\begin{equation}\label{eqn:Elman2}
\max_{v\in {\mathbb R}^n \atop \|v\|=1} \min_{p\in \piF{R}} {\|p(A)v\|} \leq
\left(1-\frac{\nu(A)^2}{\|A\|^2}\right)^{k/2}.
\end{equation}
It can be easily shown (see \cite{BeGoTy2005}), that the bound \eqref{eqn:Elman2} holds for general nonsingular matrices {$A\in {\mathbb C}^{n\times n}$},
without any assumption on the Hermitian part of $A$.

Starke proved in~\cite[Section~2.2]{T:St1994} and the subsequent paper~\cite[Theorem~3.2]{St1997},
that if~$A\in {\mathbb R}^{n\times n}$ has a positive definite symmetric part $M$, then
\begin{equation}\label{eqn:Starke}
\max_{v\in {\mathbb R}^n \atop\|v\|=1} \min_{p\in \piF{R}} \|p(A)v\|
\leq
\left(1-\nu(A)\nu(A^{-1})\right)^{k/2}.
\end{equation}
%
%If $A$ is real and $M$ is positive definite, then $\nu(A)=\lambda_{\min}(M)$, and
For a general nonsingular matrix we have
\begin{align}\label{eq:mu}
%\frac{\lambda_{\min}(M)}{\lambda_{\max}(A^TA)}
{\frac{\nu(A)}{\|A\|^2}}
&\leq \min_{{w\in\mathbb{C}^n\setminus\{0\}}}
\left|\frac{\sxy{Aw}{w}}{\sxy{w}{w}}\,\frac{\sxy{w}{w}}{\sxy{Aw}{Aw}}\right|
=\min_{{v\in\mathbb{C}^n\setminus\{0\}}}\left|\frac{\sxy{A^{-1} v}{v}}{\sxy{v}{v}}\right|=
\nu(A^{-1}),
\end{align}
which yields
$$1-\nu(A)\nu(A^{-1})\leq {1-\frac{\nu(A)^2}{\|A\|^2}.}
%1-\frac{\lambda_{\min}(M)^2}{\lambda_{\max}(A^TA)}.
$$
Hence, as pointed out by Starke in~\cite{T:St1994,St1997}, the bound (\ref{eqn:Starke})
improves Elman's bound~\eqref{eqn:Elman}. In~\cite[Corollary~6.2]{EiEr2001}, Eiermann and
Ernst proved that the bound \eqref{eqn:Starke} holds for any nonsingular
matrix $A\in {\mathbb C}^{n\times n}$. In particular, no assumption on the Hermitian part
of $A$ is required.
Note, however, that the bound \eqref{eqn:Starke} provides some information about the
convergence of (worst-case) GMRES only when $0\notin\mathcal{F}(A)$, or, equivalently, $0\notin\mathcal{F}(A^{-1})$.

In many situations the convergence of GMRES and even of
worst-case GMRES is superlinear, and therefore linear bounds like \eqref{eqn:Elman2}
and \eqref{eqn:Starke} may significantly overestimate the (worst-case) GMRES residual
norms. Nevertheless, such bounds can be very useful in the practical analysis of the
GMRES convergence, since they depend only on simple properties of the matrix $A$, which
may be estimated also in complicated applications. For example, Starke used his bound
in~\cite{T:St1994,St1997} to analyze the dependence of the convergence of hierarchical
basis and multilevel preconditioned GMRES applied to finite element discretized elliptic
boundary value problems on the mesh size and the size of the skew-symmetric part of
the preconditioned discretized operator. Similarly, Elman's bound was used in
the analysis of the GMRES convergence for finite element discretized elliptic
boundary value problems that are preconditioned with additive and multiplicative
Schwarz methods~\cite{CaWi1992,CaWi1993}. Many further such applications exist.

A straightforward upper bound on the $k$th worst-case GMRES residual norm is given by the
$k$th {\em ideal GMRES approximation}, originally introduced in~\cite{GrTr1994},
\begin{equation}\label{eqn:GMRESbounds}
\underbrace{\max_{v\in {\mathbb C}^n \atop \|v\|=1}
\min_{p\in \piF{F}}{\|p(A)v\|}}_{\mbox{\small worst-case GMRES}}
\leq
\min_{p\in \piF{F}}\max_{v\in {\mathbb C}^n \atop \|v\|=1} \|p(A)v\|=
\underbrace{\min_{p\in \piF{F}} \|p(A)\|}_{\mbox{\small ideal GMRES}}.
\end{equation}
As shown by examples in~\cite{FaJoKnMa1996,To1997} and more recently in~\cite{FaLiTi2013},
there exist matrices $A$ and iteration steps $k$ for which the inequality in~(\ref{eqn:GMRESbounds})
can be strict. The example in~\cite{To1997} even shows that the ratio of worst-case
and ideal GMRES can be arbitrarily small. A survey of the mathematical relations between
the two approximation problems in \eqref{eqn:GMRESbounds} is given
in the introductory sections of~\cite{TiLiFa2007}.

The main goal in this paper is to show that the right hand side of the bound \eqref{eqn:Starke}
also represents an upper bound on the ideal GMRES approximation for general (nonsingular)
complex matrices. This has been stated without proof {already}
in our paper~\cite[p.~168]{LiTi2005} and later
in the book~\cite[Section~5.7.3]{B:LiSt2013}. In light of the practical relevance
of Elman's and Starke's bounds, and of the fact that the inequality in~(\ref{eqn:GMRESbounds})
can be strict, we believe that providing a complete proof is important. This proof and a
further discussion of the bounds are given in Section~\ref{sec:proof}. In Section~\ref{sec:improve} we derive some improvements of the {considered} bounds.

Throughout the rest of this paper we will consider the general setting
\cblue{in ${\mathbb C}$}.

\section{Proof of the ideal GMRES bound}\label{sec:proof}

Consider a nonsingular matrix $A\in {\mathbb C}^{n\times n}$, a unit norm vector $v\in {\mathbb C}^n$,
and the minimization problem
\[
 \min_{\alpha\in\mathbb{C}}\|v-\alpha Av \|^{2}.
\]
It is easy to show that the minimum is attained for
\[
\alpha_{*}\equiv \frac{{\langle v, Av\rangle}}{\langle Av,Av\rangle},
\]
and that
\begin{eqnarray}\label{eqn:gmres1}
 \|v-\alpha_{*} Av \|^{2}&=&1-\frac{\langle Av,v\rangle}{\langle v,v\rangle}
 \frac{\langle v, Av\rangle}{\langle Av,Av\rangle}.
 %= 1-\frac{{\langle A^{-1} w , w\rangle}}{\langle w,w\rangle}\frac{\langle Av,v\rangle}{\langle v,v\rangle},
\end{eqnarray}
%where $w\equiv Av$.
%
Another result we will use below is that ideal and worst-case GMRES for any matrix $A\in {\mathbb C}^{n\times n}$
are equal in the {iteration} step $k=1$, i.e.,
\begin{equation} \label{eqn:minmax}
\min_{\alpha\in\mathbb{C}}\max_{v\in {\mathbb C}^n\atop \|v\|=1} \|v-\alpha Av\|
= \max_{v\in {\mathbb C}^n\atop \|v\|=1} \min_{\alpha\in\mathbb{C}} \|v-\alpha Av\|;
\end{equation}
see~\cite[Theorem~1]{Jo1994} and~\cite[Theorem~2.5]{GrGu1994}. This equality has also been shown in the
context of bounded linear operators on a Hilbert space; see~\cite{AsPt1971} or~\cite[Section~3.2]{B:GuRaDu1997}
and the references given there.

After these preparations we can now state and prove our main result.

\medskip
\begin{theorem}\label{thm:main}
If $A\in\mathbb{C}^{n\times n}$ is nonsingular, then for all $k\geq 1$ we have
\begin{equation}\label{eqn:final}
\min_{p\in\piF{C}} \|p(A)\|\leq \left(1-\nu(A)\nu(A^{-1})\right)^{k/2}.
\end{equation}
Moreover, if $M=\frac12(A+A^H)$ is positive definite, then
\begin{equation}\label{eqn:final2}
\min_{p\in\piF{C}} \|p(A)\|\leq
\left(1-\frac{\lambda_{\min}(M)^2}{\lambda_{\max}(A^HA)}\right)^{k/2}.
\end{equation}
\end{theorem}

\begin{proof}
The ideal GMRES approximation satisfies
\begin{equation}\label{eqn:linearbound}
\min_{p\in\piF{C}} \|p(A)\| \leq
\min_{\alpha \in\mathbb{C}} \|(I-\alpha A)^k\| \leq
\min_{\alpha \in\mathbb{C}} \|I-\alpha A\|^k.
\end{equation}
Using \eqref{eqn:minmax} and \eqref{eqn:gmres1} we then get
\begin{align}
\min_{\alpha\in\mathbb{C}} \|I-\alpha A\|^k &=
\min_{\alpha\in\mathbb{C}}\max_{v\in {\mathbb C}^n\atop\|v\|=1} \|v-\alpha Av\|^k %\nonumber \\
= \max_{v\in {\mathbb C}^n\atop \|v\|=1} \min_{\alpha\in\mathbb{C}} \|v-\alpha Av\|^k\nonumber\\
&= \max_{v\in {\mathbb C}^n\atop \|v\|=1} \left(\min_{\alpha\in\mathbb{C}} \|v-\alpha Av\|^2\right)^{k/2}\nonumber\\
&=\max_{v\in {\mathbb C}^n\atop \|v\|=1}
\left(1-\frac{\langle Av,v\rangle}{\langle v,v\rangle}\frac{\langle v,Av\rangle}{\langle Av,Av\rangle}\right)^{k/2}\nonumber\\
&=\left(1-\min_{v\in {\mathbb C}^n\atop \|v\|=1}\frac{\langle Av,v\rangle}{\langle v,v\rangle}\frac{\langle v,Av\rangle}{\langle Av,Av\rangle}\right)^{k/2}\nonumber\\
&\leq \left(1-
\min_{v\in {\mathbb C}^n\atop \|v\|=1}
\left|\frac{\langle Av,v\rangle}{\langle v,v\rangle}\right|
\min_{w\in {\mathbb C}^n\setminus\{0\}}
\left|\frac{{\langle A^{-1}w,w\rangle}}{\langle w,w\rangle}\right|\right)^{k/2}
\label{eqn:notequal} \\
&=\left(1-\nu(A)\nu(A^{-1})\right)^{k/2}, \nonumber
\end{align}
which proves \eqref{eqn:final}. If $M=\frac12(A+A^H)$ is positive definite, then $\lambda_{\min}(M) \leq \nu(A)$
and
$$\frac{\lambda_{\min}(M)}{\lambda_{\max}(A^HA)} \leq  \nu(A^{-1})$$
(see \eqref{eq:mu}), and then \eqref{eqn:final} implies \eqref{eqn:final2}.
\end{proof}
\medskip

The derivation of the bound \eqref{eqn:final} involves several inequalities, which are
usually not tight; see \eqref{eqn:linearbound} and \eqref{eqn:notequal}. We therefore
can expect that the right hand side of \eqref{eqn:final} is in most cases much larger
than the left hand side.

Since $0 \leq \nu(A) \nu(A^{-1}) \leq 1$ holds for every matrix $A\in {\mathbb C}^{n\times n}$,
equality holds in \eqref{eqn:final} when the ideal GMRES approximation stagnates until
the iteration step $k$, i.e., when
\begin{equation} \label{eqn:stagnation}
\min_{p\in\piF{C}}\|p(A)\|=1.
\end{equation}
For this to happen it is necessary that $0\in {\cal F}(A) \Leftrightarrow 0\in {\cal F}(A^{-1})$,
and it is necessary and sufficient that
\begin{equation*} %\label{eqn:hull}
0\in  \{z\in\mathbb{C}: \| p(A) \| \geq |p(z)| \;
\mbox{for all complex polynomials $p$ of degree $\leq k$}\}.
\end{equation*}
More information about the relation between the {\em polynomial numerical hull} (i.e. the set stated above) and the stagnation of {ideal} GMRES can be found in~\cite{FaJoKnMa1996,Gr2002}. The (complete)
stagnation of GMRES, which implies the stagnation of worst-case and ideal GMRES has been analyzed, for example, in~\cite{LiTi2004,Me2014,ZaOlEl2003}.

We can also identify some cases when one of the inequalities in \eqref{eqn:linearbound} is an equality. First note that if the left hand side of \eqref{eqn:linearbound} is larger than
zero, then the polynomial solving this minimization problem, i.e., the $k$th ideal
GMRES polynomial, is unique; see~\cite{GrTr1994,LiTi2009}. Hence, in this case the
first inequality in \eqref{eqn:linearbound} is an equality if and only if the $k$th ideal GMRES polynomial is of the form $(1-\alpha z)^k$. One of the very rare cases where this happens
without stagnation is when $A=J_\lambda$ is an $n\times n$ Jordan block with a sufficiently
large eigenvalue $\lambda>0$ and $1\leq k<n/2$; see~\cite[Theorem~3.2]{TiLiFa2007}
for details. The $k$th ideal GMRES polynomial then is $(1-\lambda^{-1} z)^k$, and we obtain
$$\min_{p\in\piF{C}} \|p(J_\lambda)\| =
\|(I-\lambda^{-1} J_\lambda)^k\|=\lambda^{-k}=\|I-\lambda^{-1} J_\lambda\|^k.$$

In this special case also the second inequality in \eqref{eqn:linearbound} is an equality.
For a more general (sufficient) criterion for equality, recall that a matrix $X\in {\mathbb C}^{n\times n}$
is called {\em radial} when its {\em numerical radius} is equal to its $2$-norm, i.e.,
\begin{equation}\label{eqn:nr}
r(X) \equiv \max_{z\in {\cal F}(X)} |z| = \|X\|.
\end{equation}
This holds if and only if $\|X^k\|=\|X\|^k$ for all $k\geq 1$. Several other equivalent
characterizations of this property are given in~\cite[Problem 27, p.~45]{B:HoJo1994};
see also~\cite{Pt1962}.
Suppose that the matrix $I-\widetilde{\alpha} A$ is radial for some
$\widetilde{\alpha} \in \mathbb C$ that solves the minimization problem
%
%\begin{equation}\label{eqn:problem2}
$\min_{\alpha\in\mathbb{C}} \|(I-\alpha A)^k\|.$
%\end{equation}
%
Then
$$\min_{\alpha\in\mathbb{C}} \|I-\alpha A\|^k\geq \min_{\alpha\in\mathbb{C}} \|(I-\alpha A)^k\|
= \|(I-\widetilde \alpha A)^k\| = \|I-\widetilde \alpha A\|^k
\geq \min_{\alpha\in\mathbb{C}} \|I-\alpha A\|^k,$$
which shows that equality holds throughout and hence also in the second
inequality in \eqref{eqn:linearbound}.

Finally, it is clear that in most cases the inequality \eqref{eqn:notequal} will be strict: When solving
\begin{equation}\label{eqn:anti}
\min_{v\in {\mathbb C}^n \atop \|v\|=1}\frac{\langle Av,v\rangle}{\langle v,v\rangle}
\frac{\langle v,Av\rangle}{\langle Av,Av\rangle} =
\min_{v\in {\mathbb C}^n \atop \|v\|=1} \left(\frac{|\langle Av,v\rangle|}{\|Av\|}\right)^2=
\min_{v \in {\mathbb C}^n} \cos^2\angle (v,Av)
\end{equation}
we try to make the vectors $v$ and $Av$ as close as possible to orthogonal, and
hence only the angle between the vectors plays a role. On the other hand, solutions of
$$
\min_{v\in\mathbb{C}^n\atop \|v\|=1}\left|\frac{\langle Av,v\rangle}{\langle v,v\rangle}\right|
\quad \mbox{and}\quad
\min_{w\in\mathbb{C}^n\setminus\{0\}}
{\left|\frac{\langle A^{-1}w,w\rangle}{\langle w,w\rangle}\right|}
$$
depend on the cosine of the angle as well as on the length of the vectors.

\section{\cblue{Notes about improvements of Theorem~\ref{thm:main}}}\label{sec:improve}
One can think about some improvements of Theorem~\ref{thm:main} when bounding
the value \eqref{eqn:anti} differently than in \eqref{eqn:notequal}.
Note that the value
\begin{equation}\label{eqn:antieigenvalue}
\mu(A) \equiv \min_{v\in {\mathbb C}^n \atop \|v\|=1} \frac{|\langle Av,v\rangle|}{\|Av\|}
\end{equation}
used in \eqref{eqn:anti} is called the first {\em total antieigenvalue} of $A$; for more details see, e.g., \cite[p.~67]{B:GuRaDu1997}. While the eigenvectors are not turned at all when multiplied by $A$, the antieigenvectors corresponding to the first total antiegenvalue (vectors solving \eqref{eqn:antieigenvalue}) are the vectors most turned by $A$. The problem of finding $\mu(A)$
has been solved for Hermitian positive definite (HPD) matrices and for normal matrices.

Using the definition of $\mu(A)$,  and the inequalities \eqref{eqn:linearbound} and \eqref{eqn:notequal}, we get
\begin{equation}\label{eqn:simplebound}
\min_{p\in\piF{C}} \|p(A)\| \leq
\min_{\alpha \in\mathbb{C}} \|I-\alpha A\|^k
=
\left(1-\mu^2(A)\right)^{k/2}.
\end{equation}
For HPD matrices $A$, the solution of \eqref{eqn:antieigenvalue} can be found
using Kantorovich inequality; see also results by Gustafson and Rao \cite{B:GuRaDu1997}. In particular, it holds that
$$
    \mu(A) = \frac{2\sqrt{\lambda_{\min}(A)\lambda_{\max}(A)}}{\lambda_{\min}(A)+\lambda_{\max}(A)}
$$
and the corresponding antieigenvector is a linear combination of the eigenvectors corresponding to the smallest and largest eigenvalue. When substituting the value $\mu(A)$ into the bound
\eqref{eqn:simplebound} we get, for HPD matrices $A$,
$$
\min_{p\in\piF{C}} \|p(A)\| \leq \left(\frac{\lambda_{\max}(A)-\lambda_{\min}(A)}{\lambda_{\max}(A) + \lambda_{\min}(A)}\right)^k.
$$

The first total antieigenvalue $\mu(A)$ can also be determined for normal matrices, see \cite{GuSe1993}, based on the knowledge of $A$'s eigenvalues
$\lambda_j = \beta_j + \mathbf{i} \delta_j $, $j=1,\dots,n$. Assuming that $0\notin\mathcal{F}(A)$, it holds that
$$
\mu^2(A) = \min_{i\neq j}\frac{\left(\beta_i|\lambda_j|+\beta_j|\lambda_i|\right)^2+\left(\delta_i|\lambda_j|+\delta_j|\lambda_i|\right)^2}{\left(|\lambda_j|+|\lambda_i|\right)^2|\lambda_i||\lambda_j|}.
$$

%
%    one can also cite the results of Friedland 1975 "On Matrix Approximation"
%    where mu is a solution or certain algebraic problem
%
For general nonormal matrices, the explicit value of $\mu(A)$ is unknown. We are interrested in finding a tight lower bound on $\mu(A)$. The lower bounds used in Starke's and Elman's bounds are
$$
    \mu(A) \geq \sqrt{\nu(A)\nu(A^{-1})} \geq \frac{\nu(A)}{\|A\|},
$$
and there is an open question whether these lower bound can be improved in general.
Notice that for HPD $A$ it holds that
$\nu(A)=\lambda_{\min}(A)$, $r(A)=\lambda_{\max}(A)$, and
\begin{equation}\label{eqn:nothold}
    \mu^2(A) = \frac{4{\nu(A)r(A)}}{(\nu(A)+r(A))^2} >
    \frac{{\nu(A)}}{r(A)}.
\end{equation}
In the following example we will show that the inequality \eqref{eqn:nothold} does not hold for general nonnormal matrices.

%Using the numerical radius defined in \eqref{eqn:nr}, we can write
%%
%$$\nu(A^{-1}) r(A)=
%\min_{v\in {\mathbb C}^n\setminus\{0\}}\left|
%\frac{\langle v,Av\rangle}{\langle Av,Av\rangle}\right|
%\max_{v\in {\mathbb C}^n\atop \|v\|=1} |\langle Av,v\rangle|.$$
%%
%If $w\in {\mathbb C}^n$ is a unit norm vector that maximizes $|\langle Av,v\rangle|$, then
%%
%$$\min_{v\in {\mathbb C}^n\setminus\{0\}}\left|
%\frac{\langle v,Av\rangle}{\langle Av,Av\rangle}\right|
%\max_{v\in {\mathbb C}^n\atop \|v\|=1} |\langle Av,v\rangle|
%\leq \left|\frac{\langle w,Aw\rangle}{\|Aw\|}\right|\left|\frac{\langle Aw,w\rangle}{\|Aw\|}\right|
%\leq 1,$$
%%
%so that $\nu(A^{-1})\leq 1/r(A)$.
%We obtain
%%
%\begin{equation}\label{eqn:varphi}
%1-\cos\beta \leq 1-\nu(A)\nu(A^{-1}).
%\end{equation}
%%
%It is tempting to think that the quantity $1-\cos\beta$ yields an improvement
%of the bound \eqref{eqn:final}. As shown in the following example, however, this
%is not the case.

\medskip
\begin{example}
Consider the matrix
\[
A=\left[\begin{array}{cc}
\lambda & 1\\
 0 & \lambda
\end{array}\right]\in {\mathbb R}^{2\times 2}\quad\mbox{with}\quad \lambda>\frac12.
\]
We have $n=2$ and we are interested in the {iteration} step $k=1$. The set ${\cal F}(A)$ is a disk with center at
$\lambda$ and radius $\frac12$ (independent of $\lambda$), so that
\[
\nu(A)=\lambda-\frac{1}{2},\quad r(A)=\lambda+\frac{1}{2},\quad
\frac{\nu(A)}{r(A)} =
\frac{2\lambda-1}{2\lambda+1}.
\]
Using \cite[Example~2]{TiLiFa2007} for $n=2$ and $k=1$,
\[
\min_{\alpha\in\mathbb{R}}\|I-\alpha A\| %=\frac{4\lambda}{4\lambda^{2}+1}
=\frac{4\lambda}{4\lambda^2+1}\quad\Rightarrow
\quad
\mu(A) = \frac{4\lambda^2-1}{4\lambda^2+1}.
\]
For the particular value $\lambda=\frac34$ we get
\[
 \frac{\nu(A)}{r(A)} = \frac{1}{5},\qquad \mu(A) = \frac{5}{13}.
\]
So, in this special case the inequality \eqref{eqn:nothold} does not hold.
\end{example}

\medskip
For a general nonnormal matrix, let us define $\beta\in (0,\frac{\pi}{2})$ by
\begin{equation}\label{eqn:beta}
\cos\beta = \frac{\nu(A)}{r(A)}.
\end{equation}
Inspired by the results of Crouzeix \cite{Cr2007},
the idea of
Beckermann, Goreinov and Tyrtyshnikov \cite{BeGoTy2005,Be2005}
was to find a simple convex and compact
domain including $\mathcal{F}(A)$, and to
bound the ideal GMRES approximation using a constant and the maximum norm
of the optimal polynomial on the inclusion domain. In particular,
for any matrix $A\in {\mathbb C}^{n\times n}$ with $0\notin {\cal F}(A)$ we have, possibly after a suitable rotation that can be done without loss of generality, the inclusion
$$\mathcal{F}(A) \subseteq \{z:\mathrm{Re}(z)\geq r(A)\cos\beta\}\,\cap\,\{|z|\leq r(A)\},$$
representing a circular segment.
Using the Riemann conformal mapping from the exterior of the inclusion domain onto the exterior of the closed unit disk, and Faber polynomials, the
authors of \cite{BeGoTy2005,Be2005} were able to prove the following result.

\medskip
\begin{theorem}\label{thm:BeckermannI}
If $A\in {\mathbb C}^{n\times n}$ is such that $0\notin\mathcal{F}(A)$, and $\beta\in(0,\frac{\pi}{2})$
is given as in \eqref{eqn:beta}, then for all $k\geq 1 $ we have
\begin{eqnarray} \label{eqn:improved}
\min_{p\in\piF{C}}\|p(A)\|
&\leq& (2+\rho_{\beta})\rho_{\beta}^{k},
\end{eqnarray}
where $\rho_{\beta} \equiv 2 \sin\left(\frac{\beta}{4-2\beta/\pi}\right)<\sin\beta$.
\end{theorem}

\medskip
{The technique introduced in \cite{BeGoTy2005,Be2005} can be used to get tighter bounds for other convex and compact inclusion domains like disks or ellipses, for which the associated Faber polynomials $F_k$ are known. Then, one can exploit the inequality
$$
\min_{p\in\piF{C}}\|p(A)\| \leq \frac{2}{|F_k(0)|}
$$
proven in \cite[(11)]{Be2005}. For example, for
a disk $D$ with center $c$ and radius $\delta$, $F_k$ is a shifted and rescaled monomial, $F_k(z)=\left(\frac{z-c}{\delta}\right)^k$, see, e.g., \cite{Cu1971}, giving the bound
\begin{equation}\label{eqn:disk0}
\min_{p\in\piF{C}}\|p(A)\|  \leq
 2 \left|\frac{\delta}{c}\right|^{k}.
\end{equation}

We will present an
alternative proof of \eqref{eqn:disk0}.
%the bound \eqref{eqn:disk0} can be proven using \cite{BaCrDe2006}.
%
It has been shown in \cite{BaCrDe2006} that if
$\mathcal{F}(A)\subseteq D$, then
\begin{equation}\label{eqn:disk1}
\|p(A)\|\leq 2\max_{z\in D}|p(z)|
\end{equation}
holds for any polynomial $p$.
Moreover, it is well known that the problem
\begin{equation*}%\label{eqn:disk}
\min_{p\in\piF{C}}\max_{z\in D}|p(z)|
\end{equation*}
is solved by the polynomial $(1-\frac{1}{c}z)^k$.
Therefore,
%
%Moreover, it is well known that the problem
%on the right-hand side of \eqref{eqn:disk1}
%
%\begin{equation*}%\label{eqn:disk}
%\min_{p\in\piF{C}}\max_{z\in D}|p(z)|
%\end{equation*}
%
%is solved by the polynomial $(1-\frac{z}{c})^k$.
%Then
\begin{equation}\label{eqn:disk}
\min_{p\in\piF{C}}\|p(A)\|  \leq
 2 \min_{p\in\piF{C}}\max_{z\in D}|p(z)| =
2\max_{z\in D}\left|1-\frac{z}{c}\right|^{k}
  = 2 \left|\frac{\delta}{c}\right|^{k}.
\end{equation}

For example, if
$$%\begin{equation}\label{eqn:params}
|c| = \frac{\nu(A)+r(A)}{2},\quad\mbox{and}\quad \delta = \frac{r(A)-\nu(A)}{2},
$$%\end{equation}
then
\begin{equation}\label{eqn:cf}
\left|\frac{\delta}{c}\right|
% = \frac{r(A)-\nu(A)}{r(A)+\nu(A)}
 =\frac{1-\cos\beta}{1+\cos\beta}
 <1-\cos\beta<
 1-\nu(A)\nu(A^{-1}).
\end{equation}
The convergence factor \eqref{eqn:cf}
reminds of the error bound for the classical Richardson iteration or the steepest descent
method; see, e.g.,~\cite[Section~5.5.2]{B:LiSt2013}. In particular, if $A$ is {Hermitian}
positive definite, then $\cos\beta=\lambda_{\min}(A)/\lambda_{\max}(A)=1/\kappa(A)$.

Also note that for any $\beta\in (0,\frac{\pi}{2})$ we have
$$\frac{1-\cos(\beta)}{1+\cos(\beta)}<
 2\sin\left(\frac{\beta}{4-\frac{2}{\pi}\beta}\right),$$
which can be verified using a mathematical software, or by a more detailed analysis.
Consequently, the convergence factor \eqref{eqn:cf} is smaller than the
convergence factor in Theorem~\ref{thm:BeckermannI}, which however is valid whenever
$0\notin {\cal F}(A)$.}

\begin{figure}
\centering
\includegraphics[width=0.7\linewidth]{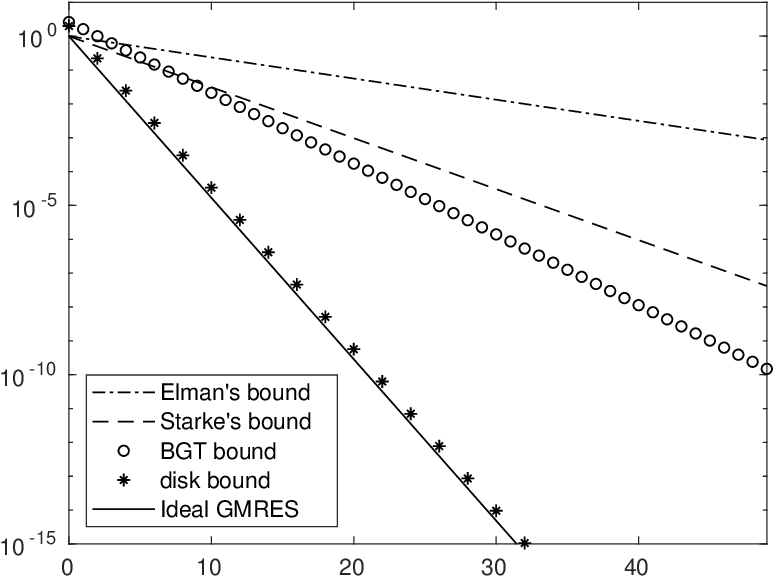}
\caption{{Ideal GMRES} and the bounds considered in this paper for a $100\times 100$ Jordan block with the eigenvalue $3$.}\label{fig:bounds}
\end{figure}

For a numerical illustration of the bounds considered in this paper we use a
single Jordan block $J_\lambda$ of the size $n=100$ and with the
eigenvalue $\lambda=3$. In Figure~\ref{fig:bounds} we plot for the first 49
iterations %the relative GMRES residual norm for a random $b$ and $x_0=0$,
{the value of the ideal GMRES approximation} (known to be $\lambda^{-k}$ in this case),
Elman's bound~\eqref{eqn:Elman}, Starke's bound~\eqref{eqn:Starke}, \cblue{the
Beckermann-Goreinov-Tyrtyshnikov} bound~\eqref{eqn:improved},
and the disk bound, i.e., the first expression in~\eqref{eqn:disk}.
We observe that the convergence factors which
determine the bounds can be quite different from each other, even in this simple
case. The disk bound is by far the best, which is due to the fact that
$\mathcal{F}(J_\lambda)$ actually is a disk centered at $\lambda$.

%In general we can of course not expect that the bounds are tight. In particular,
%the bounds are always linear while convergence of GMRES or even ideal GMRES
%is often superlinear. Nevertheless, the bounds based the field of values can
%provide a useful information in some cases.

\medskip
{\bf Acknowledgments.} This work was motivated by a question of Otto Strnad, a student of
Zden{\v e}k Strako{\v s} at the Charles University in Prague. We thank Andreas Frommer 
%\cblue{and anonymous referees
for helpful comments on previous versions of this manuscript.

%\newpage

\end{document}